\documentclass[a4paper,11pt]{article}

\usepackage[utf8]{inputenc} % allow utf-8 input
\usepackage[T1]{fontenc}    % use 8-bit T1 fonts
\usepackage{lmodern} 
\usepackage{hyperref}       % hyperlinks
\usepackage{url}            % simple URL typesetting
\usepackage{booktabs}       % professional-quality tables
\usepackage{amsfonts}       % blackboard math symbols
\usepackage{nicefrac}       % compact symbols for 1/2, etc.
\usepackage{microtype}      % microtypography

\usepackage[english]{babel}
\usepackage[utf8]{inputenc}
\usepackage{amssymb,amscd,amsfonts,bbm,mathrsfs,enumerate}
\usepackage[shortlabels]{enumitem}
\usepackage{color}
\usepackage{yhmath}
\usepackage{algorithmicx,algorithm,algpseudocode}
\usepackage{csquotes}
\usepackage{authblk}

\newtheorem{theorem}{Theorem}[section]

\newtheorem{example}{Example}
\newtheorem{remark}{Remark}

\usepackage{amsopn}
\DeclareMathOperator{\dom}{dom}

\DeclareMathOperator{\relint}{ri}
\DeclareMathOperator{\epi}{epi}
\DeclareMathOperator{\graph}{gr}
\DeclareMathOperator{\prox}{prox}

\DeclareMathOperator{\dist}{dist}

\DeclareMathOperator{\cl}{cl}

\newcommand{\eqdef}{:=}

\newcommand{\sF}{{\mathsf F}}
\newcommand{\sH}{{\mathsf H}}

\newcommand{\bP}{{{\mathbb P}}} 
\newcommand{\bE}{{{\mathbb E}}} 
\newcommand{\bN}{{{\mathbb N}}} 

\newcommand{\mC}{{\mathcal C}}

\newcommand{\sfc}{{\mathsf c}}
\newcommand{\sA}{{\mathsf A}}

\newcommand{\sB}{{\mathsf B}}
\newcommand{\sP}{{\mathsf P}}
\newcommand{\sL}{{\mathsf L}}

\newcommand{\sM}{{\mathsf M}}
\newcommand{\sG}{{\mathsf G}}

\newcommand{\maxmon}{{\mathscr M}} 
\newcommand{\Selec}{{\mathfrak S}} 

\newcommand{\mcB}{{\mathscr B}} 
 
\newcommand{\mcT}{{\mathscr T}} 
 
\newcommand{\mcF}{{\mathscr F}} 
\newcommand{\mcG}{{\mathscr G}}
\newcommand{\mcL}{{\mathscr L}}

\newcommand{\cS}{{{\mathcal S}}}

\newcommand{\cL}{{{\mathcal L}}}

\newcommand{\bR}{{{\mathbb R}}} 
\newcommand{\Hx}{{\mathcal X}}      % Space of the x variable 
\newcommand{\Hl}{{\mathcal V}}      % Space of the Lagrange multiplier lambda 
\newcommand{\Hy}{{\mathcal Y}}      % Space of the Lagrange multiplier lambda 

\newcommand{\ps}[1]{\langle #1 \rangle}
\newcommand{\wn}{{{\widetilde \nabla}}}
\usepackage[textwidth=2cm, textsize=footnotesize]{todonotes}  
\author[1]{Pascal Bianchi} 
\author[2]{Walid Hachem} 
\author[3]{Adil Salim} 
\affil[1]{LTCI, Télécom Paris, IP Paris, 75013, Paris, France.}    
\affil[2]{LIGM, CNRS, Univ.~Gustave Eiffel, F-77454 
Marne-la-Vallée, France} 
\affil[3]{Visual Computing Center, KAUST, Saudi Arabia.} 

\begin{document}

\title{A Fully Stochastic Primal-Dual Algorithm 
% \thanks{This work was partially supported by the Labex Digiteo-DigiCosme
% (OPALE project), Universit\'e Paris-Saclay.}
}


\date{27 January 2020}

\maketitle

\begin{abstract} A new stochastic primal-dual algorithm for solving a composite
optimization problem is proposed. It is assumed that all the functions /
operators that enter the optimization problem are given as statistical
expectations.  These expectations are unknown but revealed across time through
i.i.d realizations.  The proposed algorithm is proven to converge to a saddle
point of the Lagrangian function. In the framework of the monotone operator
theory, the convergence proof relies on recent results on the stochastic
Forward Backward algorithm involving random monotone operators.  An example of
convex optimization under stochastic linear constraints is considered. 
\end{abstract}

\section{Introduction} 
\label{sec:intro} 
Many applications in machine learning, statistics or signal processing require
the solution of the following optimization problem. Given two Euclidean spaces
$\Hx$ and $\Hl$, solve
\begin{equation} 
\label{pb} 
\min_{x \in \Hx} \sF(x) + \sG(x) + \sH(\sL x)
\end{equation} 
where $\sF,\sG$ and $\sH$ are lower semicontinuous convex functions such that
$\sF(x)< \infty$ for every $x$ and $\sL$ belongs to the set $\cL(\Hx, \Hl)$ of
$\Hx\to\Hl$ linear operators. 

Assuming the truth of the qualification condition $0 \in \relint ( \dom \sH -
\sL \dom \sG )$, where $\dom$ is the domain of a function and $\relint$ is the
relative interior of a set, primal-dual methods generate a sequence of primal
estimates $(x_n)_{n \in \bN}$ and a sequence of dual estimates $(\lambda_n)_{n
\in \bN}$ jointly converging to a saddle point of the Lagrangian function
$(x,\lambda) \mapsto \sF(x) + \sG(x) - \sH^\star(\lambda) + \ps{\sL x,
\lambda}$, where $\sH^\star$ is the Fenchel conjugate of $\sH$. There is a rich
literature on such algorithms which cannot be exhaustively
listed~\cite{chambolle2011first,vu2013splitting,con-jota13}. 

In this paper, it is assumed that the quantities that enter the minimization
problem are unavailable or difficult to compute numerically, and have to be
replaced with random quantities.  Specifically, let $(\Xi,\mcG,\mu)$ be a
probability space, and let $f : \Xi \times \Hx \rightarrow \bR$ and $g : \Xi
\times \Hl \rightarrow (-\infty, +\infty]$ be two convex normal integrands (see
below).  Assume that $\sF(x)=\bE_\mu(f(\cdot,x))$ and $\sG(x) =
\bE_\mu(g(\cdot,x))$. In addition, let $L$ be a measurable function from
$(\Xi,\mcG,\mu)$ to $\cL(\Hx,\Hl)$ (\textit{i.e} a random matrix), and assume
that $\sL = \bE_\mu L(\cdot)$. Finally, assume that $\sH^\star$ takes the form
$\sH^\star(\lambda) = \bE_\mu (p(\cdot,\lambda))$, where $p$ is a normal convex
integrand.  In order to solve Problem~\eqref{pb}, no one of the objects $\sF$,
$\sG$, $\sH$ and $\sL$ is available. Instead, the observer is given the
functions $f$, $g$, $p$, and $L$, along with a sequence of independent and
identically distributed (i.i.d.) random variables $(\xi_n)$ with the
probability distribution $\mu$. In this paper, a new stochastic primal dual
algorithm based on this data is proposed to solve this problem.  The
convergence proof for this algorithm relies on the monotone operator theory.
The algorithm is built around an instantiation of the stochastic
Forward-Backward (FB) algorithm involving random monotone operators that was
introduced in~\cite{bia-hac-16}. It is proven that the weighted means of the
iterates of the algorithm, where the weights are given by the step sizes of the
algorithm, converges almost surely to a saddle point of the Lagrangian
function.

To our knowledge, the proposed algorithm is the first method that
allows to solve Problem~\eqref{pb} in a fully stochastic setting with weak
assumptions on the noise.  Existing methods typically allow to handle
subproblems of Problem~\eqref{pb} in which some quantities used in this problem
are assumed to be available or set to
zero~\cite{ouyang2013stochastic,rosasco2015stochastic,toulis2015stable,yu2017online}.
In particular, the new algorithm generalizes the stochastic gradient algorithm,
the stochastic proximal point
algorithm~\cite{patrascu2017nonasymptotic,toulis2015stable,bia-16},  and the
stochastic proximal gradient
algorithm~\cite{atchade2017perturbed,bia-hac-sal-(sub)jca17}. 
A close paper to ours is \cite{com-pes-16}, which deals with a FB algorithm 
with deterministic monotone operators and random additive errors. In this
reference, the convergence of the iterates is established under stringent 
summability conditions on these errors. Random block coordinate iterations 
combined with the FB algorithm were also considered in 
\cite{com-pes-15,bia-hac-iut-16,com-pes-19}.

The next section is devoted to rigorously stating the problem and the main
result. An application example is also considered.  Section~\ref{sec:proof} is
devoted to the proof of our main theorem.

\paragraph*{Some notations.} 

The notation $\mcB(\Hx)$ will refer to the Borel $\sigma$-field of $\Hx$.  Both
the operator norm and the Euclidean vector norm will be denoted as $\|\cdot\|$.
The distance of a point $x$ to a set $S$ is denoted as $\dist(x, S)$. As
mentioned above, we denote as $\cL(\Hx, \Hl)$ the set of linear operators,
identified with matrices, from $\Hx$ to $\Hl$.  The set of proper, lower
semicontinuous convex functions on $\Hx$ is $\Gamma_0(\Hx)$. The set of
real-valued $k$--summable sequences is $\ell^k$. 

\section{Problem description and main result}
\label{sec:pb} 

We start by recalling some mathematical definitions.  Let $(\Xi, \mcG,
\mu)$ be a probability space where the $\sigma$-field $\mcG$ is $\mu$-complete,
and let $\Hx$ be an Euclidean space.  A function $h: \Xi \times \Hx \to
(-\infty, \infty]$ is said a convex normal integrand~\cite{roc-wets-livre98} if
$h(s,\cdot)$ is convex, and if the set-valued mapping $s\mapsto \epi h(s,
\cdot)$ is closed-valued and measurable in the sense
of~\cite[Chap.~14]{roc-wets-livre98}, where $\epi$ is the epigraph of a
function. We shall always assume that $h(s,\cdot) \in \Gamma_0(\Hx)$ for 
$\mu$--almost all $s \in \Xi$. 
Given $x \in \Hx$, denote as $\partial h(s,x)$ the subdifferential of 
$h(s,\cdot)$ at $x$. For $r \in [1, \infty)$, let $\mcL^r(\mu)$ be the space 
of the ${\mcG}$-measurable functions $\varphi : \Xi \to \Hx$ such that 
$\int \| \varphi \|^r d\mu < \infty$.  If 
$\mu(\{s \in \Xi \, : \, \partial h(s,x) \neq \emptyset \}) < 1$, 
set $\Selec^r_{\partial h(\cdot,x)} \eqdef \emptyset$, otherwise, 
\[
\Selec^r_{\partial h(\cdot,x)} \eqdef 
\{ \varphi \in \mcL^r(\mu) \, : \, 
\varphi(s) \in \partial h(s,x) \ \mu-\text{almost everywhere (a.e.)} \} 
\]
is the set of the so-called $r$--integrable selections of the measurable 
set-valued function $s\mapsto \partial h(s, x)$. Denoting as $\cl$ the closure
of a set, the so-called selection integral of $\partial h(\cdot,x)$ is the set
\begin{equation}
\label{selint} 
\bE_\mu \partial h(\cdot,x) \eqdef \cl{\left\{ \int_\Xi \varphi d\mu \ : \ 
  \varphi \in \Selec^1_{\partial h(\cdot,x)} \right\}} 
\end{equation}
that might be empty. 
Note that we use the same notation $\bE_\mu$ for these set-valued expectations
and for the classical single-valued expectations. 
% The meaning of $\bE_\mu$ will be always clear from the context. 

We now state our problem.  Let $f: \Xi \times \Hx \to (-\infty, \infty]$ be a
convex normal integrand, assume that $\bE_\mu | f(\cdot,x) | < \infty$ for all
$x \in \Hx$, and consider the convex function $\sF(x) \eqdef \bE_\mu f(\cdot,
x)$ which domain is $\Hx$.  Let $g: \Xi \times \Hx \to (-\infty, \infty]$ be a
convex normal integrand, and let $\sG(x) \eqdef \bE_\mu g(\cdot, x)$, where
the integral $\bE_\mu$ is defined as the sum 
\[
\int_{\{s \, : \, g(s, x) \in [0,\infty) \}} 
           g(s, x) \, \mu(ds)  + 
\int_{\{s \, : \, g(s, x) \in ]-\infty, 0[\}} 
           g(s, x) \, \mu(ds)  +  I(x) \, , 
\]
and 
\[
I(x) = \left\{\begin{array}{cl} + \infty, &\text{if } 
         \mu(\{s : g(s, x) = \infty \}) > 0, \\
0, &\text{otherwise} \, , \end{array}\right.  
\]
and where the convention $(+\infty) + (-\infty) = + \infty$ is used. The
function $\sG$ is a lower semi continuous convex function if $\sG(x) > -\infty$
for all $x$, which we assume. We shall assume that $\sG$ is proper.  In a
similar manner, let $p : \Xi  \times \Hl \to (-\infty,\infty]$ be a convex
normal integrand, assume that $\sP :\lambda\mapsto \bE_\mu p(\cdot,\lambda)$ belongs to
$\Gamma_0(\Hl)$, and let $\sH$ be its Fenchel conjugate 
(thus, $\sH^\star = \sP$). Finally, let $L : \Xi \to \cL(\Hx,\Hl)$ be an 
operator-valued measurable function, assume that $\| L \|$ is $\mu$-integrable, and let $\sL \eqdef \bE_\mu L$.

Having introduced these functions, our purpose is to find a solution 
$x \in \Hx$ of Problem~\eqref{pb},  
where the set of such points is assumed non empty. 
To solve this problem, the observer is given the functions 
$f, g, p, L$, and a sequence of i.i.d random variables 
$(\xi_n)_{n\in\bN}$ from a probability space $(\Omega, \mcF, \bP)$ to 
$(\Xi, \mcG)$ with the probability distribution $\mu$. 

Denote as $\prox_{h}(x) \eqdef \arg\min_{y \in \Hx} h(y) + \|y-x\|^2 / 2$ the
Moreau's proximity operator of a function $h \in \Gamma_0(\Hx)$.  We also
denote as $\partial_0 h(x)$ the least norm element of the set 
$\partial h(x)$, which is known to exist and to be 
unique \cite{bau-com-livre11}. 
% As is well known, 
% $\partial_0 h(x) = \lim_{\gamma\downarrow 0} \nabla h_\gamma(x)$, where
% $h_\gamma$ is Moreau's envelope of $h$ for $\gamma > 0$, which is the
% differentiable function defined as  
% $h_\gamma(x) = \min_{y \in \Hx} h(y) + \|y-x\|^2 / (2\gamma)$. 
Similarly, $\partial_0 f(s,x)$ will refer to the least
norm element of $\partial f(s,x)$ which was introduced above. We shall also
denote as $\wn f(s,x)$ a measurable subgradient of $f(s,\cdot)$ at $x$. 
Specifically, 
$\wn f : (\Xi \times \Hx, \mcG \otimes \mcB(\Hx)) \to (\Hx, \mcB(\Hx))$ is a 
measurable function such that for each $x\in \Hx$, 
$\wn f(\cdot, x) \in \Selec^1_{\partial f(\cdot,x)}$, which is known to be 
non empty thanks to the integrability assumption 
$\bE_\mu | f(\cdot,x) | < \infty$ \cite{roc-wet-82}. 
A possible choice for $\wn f(s, x)$ is $\partial_0 f(s,x)$ 
\cite[\S2.3 and \S3.1]{bia-hac-16}.  
Turning back to Problem~\eqref{pb}, our purpose will be to find a saddle point
of the Lagrangian 
$(x,\lambda) \mapsto \sF(x) + \sG(x) - \sH^\star(\lambda) + \ps{\sL x, \lambda}$. Denoting as $\cS \subset \Hx
\times \Hl$ the set of these saddle points, an element $(x, \lambda)$ of
$\cS$ is characterized by the inclusions 
\begin{equation}
\label{selle}
  \left\{
    \begin{array}[h]{lccl}
      0 & \in & \partial \sF(x) + \partial \sG(x)& + \sL^T\lambda , \\
      0 & \in & -\sL x& + \partial \sH^\star(\lambda)\,. 
    \end{array}
\right.
% \nonumber 
\end{equation}

Consider a sequence of positive weights $(\gamma_n)_{n\in\bN}$. The algorithm
proposed here consists in the following iterations applied to the random vector
$(x_n,\lambda_n) \in \Hx \times \Hl$. 

\begin{equation} 
\label{theta=0}
  \begin{split}
    x_{n+1} &= \prox_{\gamma_{n+1} g(\xi_{n+1},\cdot)}
  \left(
    x_n -\gamma_{n+1}(\wn f(\xi_{n+1}, x_n)+L(\xi_{n+1})^T \lambda_n)
     \right),  \\
    \lambda_{n+1} &= \prox_{\gamma_{n+1} p(\xi_{n+1},\cdot)} \left( \lambda_n +
    \gamma_{n+1}
    L(\xi_{n+1}) x_{n} \right)  \,.
  \end{split}
\end{equation}

The convergence of Algorithm~\eqref{theta=0} is stated by the next theorem in terms of weighted averaged estimates 
\[
\bar x_n = \frac{\sum_{k=1}^n \gamma_k x_k}{\sum_{k=1}^n \gamma_k} , \ 
\text{and} \ 
\bar \lambda_n = 
  \frac{\sum_{k=1}^n \gamma_k \lambda_k}{\sum_{k=1}^n \gamma_k}. 
\]

\begin{theorem}
\label{th:theta0}
Consider Problem~\eqref{pb}, and let the following assumptions hold. 

\begin{enumerate} 
\item\label{step} The step size sequence satisfies 
 $(\gamma_n) \in \ell^2 \setminus \ell^1$, and 
 $\gamma_{n+1} / \gamma_n \rightarrow 1$ as $n\to\infty$. 
% \item\label{empty} The set $\cS$ is not empty.\asnote{Intersection avec l'hyp qui suit}

\item\label{interchange} The function $\sG$ satisfies 
$\partial \sG(x) = \bE_\mu \partial g(\cdot,x)$ for each $x \in \Hx$. 

\item\label{sm} There exists an integer $m\geq 2$ that satisfies the following 
 conditions: 
 \begin{itemize} 
  \item The function $L$ is in $\mcL^{2m}(\mu)$. 
  \item There exists a point $(x_\star, \lambda_\star) \in \cS$, and 
  three functions 
  $\varphi_f \in \Selec^{2m}_{\partial f(\cdot, x_\star)}$,
  $\varphi_g \in \Selec^{2m}_{\partial g(\cdot, x_\star)}$, and    
  $\varphi_p \in \Selec^{2m}_{\partial p(\cdot, \lambda_\star)}$ such that 
  \begin{equation}
  \label{phif} 
  \bE_\mu \varphi_f +
  \bE_\mu \varphi_g + \sL^T \lambda_\star = 0, \ \text{and } \ -\sL x_\star +
  \bE_\mu \varphi_p = 0. 
  \end{equation}
 \end{itemize} 
 Moreover, for every point $(x_\star, \lambda_\star) \in \cS$, there exist 
  three functions 
  $\varphi_f \in \Selec^{2}_{\partial f(\cdot, x_\star)}$,
  $\varphi_g \in \Selec^{2}_{\partial g(\cdot, x_\star)}$, and    
  $\varphi_p \in \Selec^{2}_{\partial p(\cdot, \lambda_\star)}$ such that~\eqref{phif} holds.

\item\label{d0} For any compact set $K \subset\dom\partial \sG$, there exist
 $\varepsilon \in (0,1]$ and $x_0 \in \dom\partial \sG$ such that 
\begin{equation*}
\sup_{x \in K} \bE_\mu \|\partial_0 g(\cdot, x)\|^{1+\varepsilon} < +\infty, 
\ \text{and } \ 
\bE_\mu \|\partial_0 g(\cdot, x_0)\|^{1+1/\varepsilon} < +\infty.
\end{equation*}

\item\label{croissance} There exists a measurable $\Xi \to \bR_+$ function 
$\beta$ such that $\beta^{2m}$ is $\mu$-integrable, where $m$ is the integer
provided by Assumption~\ref{sm}, and such that for all $x \in \Hx$,
$$
\|\wn f(s, x)\| \leq \beta(s) (1+\|x\|).
$$
Moreover, there exists a constant $C > 0$ such that 
$\bE_\mu \|\wn f(\cdot, x)\|^4 \leq C (1+\|x\|^{2m})$. \\

\item\label{reglin} Writing $D_{\partial g}(s) = \dom \partial g(s, \cdot)$, 
there exists $C >0$ such that for all $x \in \Hx$,
\[
\bE_\mu \dist( x, D_{\partial g}(\cdot))^2 \geq C \dist(x, \dom \partial \sG)^2. 
\]

\item\label{projo} There exists $C > 0$ such that for any $x \in \Hx$ and any 
$\gamma >0$,
\[
\int \|\prox_{\gamma g(s, \cdot)}(x) - \Pi_{g}(s,x)\|^4 \mu(ds) 
\leq C \gamma^4 (1+\|x\|^{2m}), 
\]
where $\Pi_g(s,\cdot)$ is the projection operator onto
$\cl(\dom \partial g(s, \cdot))$, and where $m$ is the integer provided by
Assumption~\ref{sm}. 
\item\label{p==g}
Assumptions \ref{interchange}, \ref{d0}, \ref{reglin} and \ref{projo} 
hold true when the function $g$ is replaced with $p$ and the space $\Hx$ is replaced with $\Hl$. 

\end{enumerate}

Then, the sequence $(x_n,\lambda_n)$ is bounded in $\mcL^{2m}(\Omega)$ and the sequence $(\bar x_n, \bar \lambda_n)$ converges almost 
surely (a.s.) to a random variable $(X,\Lambda)$ supported by $\cS$. 
\end{theorem}

Let us now discuss our assumptions. Assumption~\ref{step} is standard in 
the decreasing step case. 
Assumption~\ref{interchange} requires that the interchange of the 
expectation $\bE_\mu g(\cdot, x)$ and the subdifferentiation be possible. 
Let us provide some sufficient conditions for this to be true. 
By \cite{roc-wet-82}, this will be the case if the following conditions hold: 
\emph{i)} the set-valued mapping $s\mapsto \cl\dom g(s, \cdot)$ is constant 
$\mu$-a.e., where $\dom g(s,\cdot)$ is the domain of $g(s,\cdot)$, 
\emph{ii)} $\sG(x) < \infty$ whenever $x \in \dom g(s,\cdot)$ $\mu$-a.e., 
\emph{iii)} there exists $x_0 \in \Hx$ at which $\sG$ is finite and continuous. 
Another case of practical importance where this interchange is permitted is 
the following. Let $m$ be a positive integer, and let $\mC_1, \ldots \mC_m$ be 
a collection of closed and convex subsets of $\Hx$. Let 
$\mC = \cap_{i=k}^m \mC_k$ be non empty, and assume
that the normal cone $N_{\mC}(x)$ of $\mC$ at $x$ satisfies the identity 
$N_{\mC}(x) = \sum_{k=1}^m N_{\mC_k}(x)$ for each $x\in \Hx$, where the 
summation is the usual set summation. As is well known, this identity holds 
true under a qualification condition of the type 
$\cap_{k=1}^m \relint \mC_k \neq \emptyset$ (see also 
\cite{bauschke1999strong} for other conditions). Now, assume that   
$\Xi = \{1,\ldots, m\}$ and that $\mu$ is an arbitrary probability measure
putting a positive weight on each $\{k\} \subset \Xi$. Let $g(s,x)$ be the 
indicator function 
\begin{equation} 
\label{intercpct} 
g(s,x) = \iota_{\mC_s}(x) \ \text{for} \  (s,x) \in \Xi \times \Hx.
\end{equation}  
Then it is obvious that $g$ is a convex normal integrand, $\sG = \iota_{\mC}$, 
and $\partial \sG(x) = \bE_\mu \partial g(\cdot,x)$. 
We can also combine these two types of conditions: let $(\Sigma, \mcT, \nu)$ 
be a probability space, where $\mcT$ is $\nu$-complete, and let 
$h:\Sigma \times \Hx \to (-\infty, \infty]$ be a
convex normal integrand satisfying the conditions \emph{i)}--\emph{iii)} above.
Consider the closed and convex sets $\mC_1,\ldots,\mC_m$ introduced above, and
let $\alpha$ be a probability measure on the set $\{0,\ldots, m \}$ such that
$\alpha(\{k\}) > 0$ for each $k \in \{0,\ldots, m \}$. Now, set 
$\Xi = \Sigma \times \{0,\ldots,m\}$, $\mu = \nu \otimes \alpha$, and 
define $g : \Xi \times \Hx \to (-\infty, \infty]$ as
\[
g(s, x) =  \left\{\begin{array}{ll} \alpha(0)^{-1} h(u,x) &
\text{if } k = 0, \\
\iota_{\mC_k}(x) & \text{otherwise}, 
\end{array}\right. 
\]
where $s = (u, k) \in \Sigma \times \{0,\ldots,m\}$. Then it is clear that 
\[
\sG(x) = \frac{1}{\alpha(0)} \int_\Sigma h(u, x) \nu(du) + \iota_{\mC}(x) \, , 
\]
and 
\[
\partial \sG(x) = \bE_\mu \partial g(\cdot,x) = 
 \frac{1}{\alpha(0)} \bE_\nu \partial h(\cdot, x) + 
 \sum_{k=1}^m N_{\mC_k}(x) \, . 
\]
Assumption~\ref{sm} is a moment assumption that is
generally easy to check. Note that this assumption requires the set of
saddle points $\cS$ to be non empty. Notice the relation between 
Equations~\eqref{phif} and the two inclusions in~\eqref{selle}. Focusing on 
the first inclusion and using Assumption~\ref{interchange}, there exist 
$a \in \partial \sF(x_\star) = \bE_\mu\partial f(\cdot, x_\star)$ and 
$b \in \partial \sG(x_\star) = \bE_\mu\partial g(\cdot, x_\star)$ such that
$0 = a + b + \sL^T \lambda_\star$. Then, Assumption~\ref{sm} states that $a$ 
and $b$ can be taken in such a way that there are two measurable
selections $\varphi_f$ and $\varphi_g$ of 
$\partial f(\cdot, x_\star)$ and $\partial g(\cdot, x_\star)$ respectively
which are both in $\mcL^{2m}(\mu)$ and which satisfy $a = \bE_\mu\varphi_f$ and
$b = \bE_\mu\varphi_g$. A sufficient condition for the existence of the
selections satisfying Assumption~\ref{sm} is the following~\cite{bia-hac-sal-(sub)jca17}: there exists an open
neighborhood $\mathcal N_x$ of $x_\star$ and an open neighborhood 
$\mathcal N_\lambda$ of $\lambda_\star$ such that 
$\forall x \in \mathcal N_x$, 
$\int f(s, x)^{2m} \mu(ds) < \infty$ and $\int g(s, x)^{2m} \mu(ds) < \infty$, 
and $\forall \lambda \in \mathcal N_\lambda$, 
$\int p(s, x)^{2m} \mu(ds) < \infty$. 
Note also that the larger is $m$, and the weaker is 
Assumption~\ref{projo}. 

Assumption~\ref{d0} is relatively weak and easy to check. It is interesting
to compare it with Assuption~\ref{croissance}. It is indeed much weaker than
the latter, which assumes that the growth of $\wn f(s, \cdot)$ is not
faster than linear. This is due to the fact that $g$ and $p$ enter the
algorithm~\eqref{theta=0} through the proximity operator while the function
$f$ is used explicitly in this algorithm (through its (sub)gradient). This use of the functions
$f$ is reminiscent of the well-known Robbins-Monro algorithm, where a linear growth is needed to ensure the
algorithm stability. Note that Assumption~\ref{croissance} is satisfied under the more restrictive assumption that $\nabla f(s,\cdot)$ is $L$-Lipschitz continuous without any bounded gradient assumption.  

Assumption~\ref{reglin} is quite weak, and is studied \textit{e.g} in~\cite{necoara2018randomized}. This assumption is easy to illustrate in the case
where $g(s,x) = \iota_{\mC_s}(x)$ as in~\eqref{intercpct}.
Following~\cite{bauschke1999strong}, we say that the subsets
$(\mC_1,\dots,\mC_m)$ are linearly regular if there exists $C>0$ such that for
every $x$, 
\[
\max_{i=1\dots m} \dist(x,\mC_i)\geq C \dist(x,\mC) . 
\]
Sufficient conditions for a collection of sets to satisfy the above condition
can be found in \cite{bauschke1999strong} and the references therein. Note that
this condition implies that $N_{\mC}(x) = \sum_{i=1}^m N_{\mC_i}(x)$. 
Let us finally discuss Assumption~\ref{projo}. As $\gamma\to 0$, it is known
that $\prox_{\gamma g (s,\cdot)}(x)$ converges to $\Pi_g(s,x)$ for every 
$(s,x)$. Assumption~\ref{projo} provides a control on the convergence rate. 
This assumption holds under the sufficient condition that for $\mu$-almost 
every $s$ and for every $x\in \dom \partial g(s,\cdot)$,
\[
\|\partial g_0(s,x)\|\leq \beta(s)(1+\|x\|^{m/2})\, , 
\]
where $\beta$ is a positive random variable with a finite fourth 
moment~\cite{bia-16}.

We now consider an application example of Theorem~\ref{th:theta0}.

\begin{example}
Let $\sfc \in \Hl$. Setting $\sH = \iota_{\{\sfc\}}$, where $\iota_{\mC}$ is
the indicator function of the set $\mC$, Problem~\eqref{pb} boils down to the
linearly constrained problem
\begin{equation} 
\label{pbcons} 
\min_{x \in \Hx} \sF(x) + \sG(x) \quad \text{s.t.} \quad \sL x = \sfc.
\end{equation}
If we assume that $\sfc = \bE_\mu(c(\cdot))$ where 
$c(\cdot) : \Xi \rightarrow \Hl$ is a random vector, then our problem amounts
to randomizing the constraints and to handling these stochastic constraints 
online. Such a context is encountered in various fields of machine learning, 
as the Neyman-Pearson classification, or in online so-called Markowicz 
portfolio optimization.  

Since $\sH^\star(\lambda) = \ps{\lambda,\sfc}$, we simply need to put
$p(\cdot,\lambda) = \ps{\lambda,c(\cdot)}$, and Algorithm~\eqref{theta=0}
becomes: 
\begin{align*} 
    x_{n+1} &= \prox_{\gamma_{n+1} g(\xi_{n+1},\cdot)}
  \left(
    x_n -\gamma_{n+1}(\wn f(\xi_{n+1}, x_n)+L(\xi_{n+1})^T \lambda_n)
     \right),  \\
    \lambda_{n+1} &= \lambda_n +
    \gamma_{n+1}
    \left(L(\xi_{n+1}) x_{n} - c(\xi_{n+1}) \right)  \,.
\end{align*} 
To go further, let us particularize Problem~\eqref{pbcons} to the case of the
Markowicz portfolio optimization, and check the assumptions of
Theorem~\ref{th:theta0} to complete the picture. In this case, $\xi$ is a
$\Hx$--valued random variable with a second moment, $\sF(x) =
\bE_\mu\ps{x,\xi}^2$, $\sG(x) = \iota_{\Delta}(x)$ where $\Delta$ is the
probability simplex, $\sL = \bE_\mu(\xi^{T})$, and $\sfc$ is some real positive
number. Note that it is usually assumed that $\sL = \bE_\mu(\xi^{T})$ is fully
known or estimated, which we don't do here. We of course assume that the 
qualification condition $\sfc \in \relint \sL \Delta$ holds true. 

Assumptions~\ref{interchange} and~\ref{d0} of the statement of
Theorem~\ref{th:theta0} are immediate for both $g$ and $p$. One can check that
Assumption~\ref{sm} is satisfied for $m=2$ if we assume that $\bE_\mu \| \xi
\|^4 < \infty$, which also ensures the truth of Assumption~\ref{croissance}.
Assumptions~\ref{reglin} and~\ref{projo} are trivially satisfied for $g$ and
$p$, since $\prox_{\gamma g(s,\cdot)} = \Pi_g(s,\cdot)$, and since $p(s,\cdot)$
has a full domain. 

\end{example}

\section{Proof of Theorem~\ref{th:theta0}} 
\label{sec:proof}

The proof of Theorem~\ref{th:theta0} makes use of the monotone operator theory.
We begin by recalling some basic facts on monotone operators. All the results
below can be found in \cite{bre-livre73,bau-com-livre11} without further
mention.

A set-valued mapping $\sA : \Hx \rightrightarrows \Hx$ on the Euclidean space
$\Hx$ will be called herein an operator. An operator with singleton values is identified with a function. As above, the domain of $\sA$  
is $\dom(\sA) = \{ x \in \Hx \, : \, \sA(x) \neq \emptyset \}$. The graph
of $\sA$ is $\graph(\sA) = \{ (x,y) \in \Hx\times\Hx \, : \, y \in \sA(x) \}$.  
The operator $\sA$ is said monotone if 
$\forall (x,y),(x',y')\in \graph(\sA)$, 
$\ps{y-y',x-x'}\geq 0$. A monotone operator with non empty domain is said
maximal if $\graph(\sA)$ is a maximal element for the inclusion ordering in 
the family of the monotone operator graphs. 
Let $I$ be the identity operator, and let $\sA^{-1}$ be the inverse of $\sA$,
which is defined by the fact that 
$(x,y) \in \graph(\sA^{-1}) \Leftrightarrow (y,x) \in \graph(\sA)$.  An
operator $\sA$ belongs to the set $\maxmon(\Hx)$ of the maximal monotone 
operators on $\Hx$ if and only if for each $\gamma > 0$, the so-called 
resolvent $( I + \gamma \sA )^{-1}$ is a contraction defined on the whole 
space $\Hx$. In particular, it is single-valued. 
A typical element of $\maxmon(\Hx)$ is the subdifferential $\partial h$ of a
function $h \in \Gamma_0(\Hx)$. In this case, the resolvent 
$( I + \gamma\partial h)^{-1}$ for $\gamma > 0$ coincides with 
the proximity operator $\prox_{\gamma h}$. 
A skew-symmetric element of $\cL(\Hx, \Hx)$ can also be checked to be an 
element of $\maxmon(\Hx)$. 

The set of zeros of an operator $\sA$ on $\Hx$ is the set 
$Z(\sA) =  \{ x \in\Hx \, : \, 0 \in \sA(x) \}$. The sum of two operators
$\sA$ and $\sB$ is the operator $\sA + \sB$ whose image at $x$ is the
set sum of $\sA(x)$ and $\sB(x)$. 
Given two operators $\sA, \sB \in \maxmon(\Hx)$, where $\sB$ is 
single-valued with domain $\Hx$, the FB algorithm is
an iterative algorithm for finding a point in $Z(\sA + \sB)$. It reads
\[
x_{n+1} = ( I + \gamma \sA )^{-1} ( x_n - \gamma \sB(x_n) ) \, 
\]
where $\gamma$ is a positive step. 

In the sequel, we shall be interested by random elements of $\maxmon(\Hx)$ as
used in \cite{bia-16,bia-hac-16,bia-hac-sal-(sub)jca17}. A random element
of $\maxmon(\Hx)$ is a measurable function $M : \Xi \to\maxmon(\Hx)$ in the
sense of~\cite{att-79}, where $(\Xi, \mcG, \mu)$ is the probability space
introduced at the beginning of Section~\ref{sec:pb}. 
In particular, when $h:\Xi\times\Hx \to (-\infty,\infty]$ is a convex normal
integrand such as $h(s,\cdot)$ is proper $\mu$-a.e., 
$M(s) = \partial h(s,\cdot)$ is a random element of $\maxmon(\Hx)$. 
Moreover, when $M(s)$ is a skew-symmetric element of $\cL(\Hx,
\Hx)$ which is measurable in the usual sense (as a $\Xi \to \cL(\Hx, \Hx)$
function), then it is also a random element of $\maxmon(\Hx)$. 
If we fix $x \in \Hx$ and we denote as $M(s,x)$ its image by $M(s)$, then
the set-valued function $s \mapsto M(s,x)$ is measurable, and its 
(set-valued) expectation $\sM(x) = \bE_\mu M(\cdot, x)$ is defined similarly
to Equation~\eqref{selint} \cite{att-79,bia-16,bia-hac-16}. 
Note that $\sM$ is monotone but not necessarily maximal.

We now enter the proof of Theorem~\ref{th:theta0}. Let us set 
$\Hy = \Hx\times\Hl$, and endow this Euclidean space with the 
standard scalar product. By writing $(x,\lambda) \in \Hy$, it will be 
understood that $x \in \Hx$ and $\lambda\in \Hl$. 
For each $s \in \Xi$, define the set-valued operator $A(s)$ on $\Hy$ as 
the operator that takes $(x,\lambda)$ to 
\[
A(s, (x,\lambda) ) = \begin{bmatrix} \partial g(s,x) \\ \partial p(s,\lambda) \end{bmatrix} ,
\]
Fixing $s\in\Xi$, the operator $A(s, (x,\lambda))$ coincides with the 
subdifferential of the convex normal integrand 
$g(s,x) + p(s,\lambda)$ with respect to $(x,\lambda)$. Thus, $A(s)$ is a random
element of $\maxmon(\Hy)$. Let us also define the operator $B(s)$ as 
\[
B(s, (x,\lambda)) = 
\begin{bmatrix}
   \partial f(s, x)& + L(s)^T \lambda \\
    -L(s) x&  
\end{bmatrix} . 
\]
We can write $B(s) = B_1(s) + B_2(s)$, where 
\[
B_1(s, (x,\lambda)) = \begin{bmatrix} \partial f(s,x) \\
 0 \end{bmatrix}, \quad \text{and} \quad 
   B_2(s) = \begin{bmatrix} 0 & L(s)^T \\
                     -L(s) & 0 \end{bmatrix}
\]
($B_2(s)$ is a linear skew-symmetric operator written in a 
matrix form in $\Hy$). For each $s\in\Xi$, both these operators belong to $\maxmon(\Hy)$, and 
$\dom B_2(s) = \Hy$. Thus, $B(s) \in \maxmon(\Hy)$ by 
\cite[Cor.~24.4]{bau-com-livre11}. Moreover, since both $B_1$ and $B_2$ 
are measurable, $B$ is a random element of $\maxmon(\Hy)$. 

Since $f(\cdot, x)$ is Lebesgue-integrable for all $x\in \Hx$ by construction,
it is known that $\partial \sF(x) = \bE_\mu \partial f(\cdot, x)$ 
\cite{roc-wet-82}. Moreover, $\partial \sG(x) = \bE_\mu \partial g(\cdot, x)$ 
and $\partial \sH^\star(\lambda) = \bE_\mu \partial p(\cdot, \lambda)$ 
by Assumptions~\ref{interchange} and~\ref{p==g}. Thus, the operators  
$\sA((x,\lambda)) = \bE_\mu A(\cdot, (x,\lambda))$ and 
$\sB((x,\lambda)) = \bE_\mu B(\cdot,(x,\lambda))$ can be written as 
\[
\sA((x,\lambda)) = \begin{bmatrix} \partial \sG(x) \\ \partial \sH^\star(\lambda) \end{bmatrix}, \ \text{and } \ 
\sB((x,\lambda)) = \begin{bmatrix} \partial \sF(x)&  + \sL^T \lambda \\ -\sL x& \end{bmatrix} , 
\]
thus, these monotone operators are both maximal.  By
\cite[Cor.~24.4]{bau-com-livre11}, we also get that $\sA + \sB$ belong to
$\maxmon(\Hy)$. Moreover, recalling the system of inclusions \eqref{selle}, we
also obtain that $\cS = Z(\sA + \sB)$. 

Defining the function 
\[
b(s, (x,\lambda)) = 
\begin{bmatrix} 
   \wn f(s, x)& + L(s)^T \lambda \\
    -L(s) x& 
\end{bmatrix}  
\]
(obviously, $b(s, (x,\lambda)) \in B(s, (x,\lambda))$ $\mu$-a.e.), let us 
consider the following version of the FB algorithm
\begin{equation*} 
  (x_{n+1}, \lambda_{n+1})
= 
\left( I + \gamma_{n+1} A(\xi_{n+1}, \cdot) \right)^{-1} 
 \left( (x_n, \lambda_n) 
     - \gamma_{n+1} b(\xi_{n+1}, (x_n,\lambda_n))  \right) .
\end{equation*} 
On the one hand, one can easily check that this is exactly Algorithm~\eqref{theta=0}. 
On the other hand, this algorithm is an instance of the random FB 
algorithm studied in \cite{bia-hac-16}.  By checking the assumptions of
Theorem~\ref{th:theta0} one by one, one sees that the assumptions of
\cite[Th.~3.1 and Cor.~3.1]{bia-hac-16} are verified. Theorem~\ref{th:theta0}
follows. 

\begin{remark}
The convergence stated by Theorem~\ref{th:theta0} concerns the averaged
sequence $(\bar x_n, \bar\lambda_n)$. One can ask whether the sequence $(x_n,
\lambda_n)$ itself converges to $\cS$. This would happen if the operator
$\sA+\sB$ were so-called \textit{demipositive}~\cite{bia-hac-16}.  This happens
when, \emph{e.g.}, $\sF+\sG$ is strongly convex and $\sH$ is smooth (proof
omitted). Unfortunately, demipositivity of $\sA+\sB$ is not always guaranteed. 
\end{remark}

% \bibliographystyle{plain} 
% \bibliography{math}

\def\cprime{$'$} \def\cdprime{$''$} \def\cprime{$'$}

\end{document}